\documentclass[10pt]{amsart}
\usepackage{palatino, amssymb, amsfonts, latexsym, mathrsfs}
\usepackage{amssymb, amsmath, amscd}
\input xy
\xyoption{all}

\newtheorem{theorem}{Theorem}[section]
\newtheorem{lemma}[theorem]{Lemma}
\newtheorem{prop}[theorem]{Proposition}
\newtheorem{cor}[theorem]{Corollary}
\theoremstyle{definition}
\newtheorem{definition}[theorem]{Definition}

\theoremstyle{remark}
\newtheorem{remark}[theorem]{Remark}

\newcommand{\U}{\mathbf{U}}

\newcommand{\A}{\mathcal A}
\newcommand{\Al}{\mathcal A_\ell}
\newcommand{\e}{\varepsilon}
\newcommand{\Ud}{\dot{\mathbf U}}

\newcommand{\B}{\mathbf B}
\newcommand{\Bd}{\dot{\mathbf B}}

\newcommand{\ic}{(\underline{i},\underline{c})}
\newcommand{\icp}{(\underline{i}', \underline{c}')}
\newcommand{\icpp}{(\underline{i}'', \underline{c}'')}
\newcommand{\Pic}{\Phi_{(\underline i, \underline c)}}
\newcommand{\bPic}{\Psi_{(\underline i, \underline c)}}

\numberwithin{equation}{section}

\begin{document}

\title{Hall algebras and quantum Frobenius.}

\author{Kevin McGerty}
\address{Department of Mathematics, Imperial College London. }

\date{October, 2009.}

\begin{abstract}
Lusztig has constructed a Frobenius morphism for quantum groups at an $\ell$-th root of unity, which gives an integral lift of the Frobenius map on universal enveloping algebras in positive characteristic. Using the Hall algebra we give a simple construction of this map for the positive part of the quantum group attached to an arbitrary Cartan datum in the nondivisible case.
\end{abstract}

\maketitle

\begin{flushright}
\textit{To George Lusztig, with gratitude and admiration.}
\end{flushright}

\section{Introduction}

Let $\mathsf k$ be an algebraically field of charactistic $p>0$ and $G$ be an affine algebraic group over $\mathsf k$. Then as a variety over $\mathsf k$, the group $G$ can be equipped a Frobenius morphism $F \colon G \to G$, which by naturality is in fact an endomorphism of the algebraic group. The existence of this map is of fundamental importance in the study of the representation theory and geometry of $G$ (see for example \cite{St}, \cite{A}). 

Let $\mathbf U_v(\mathfrak g)$ be the quantum group attached to a symmetrizable Kac-Moody Lie algebra $\mathfrak g$. Lusztig \cite{L89}, \cite{L90} discovered that when the parameter $v$ is specialized to $\varepsilon$ an $\ell$-th root of unity, there is a homomorphism $Fr$ from the resulting algebra $\mathbf U_\varepsilon(\mathfrak g)$ to the integral form of the enveloping algebra $\mathcal U_\mathbb Z(\mathfrak g)$. This construction gives an integral lift of the Frobenius morphism: if $\ell= \text{char}(\mathsf k)$, then after base changing to $\mathsf k$  one obtains the transpose of the map $F$ on the hyperalgebra of $G$.

The existence of this map was fundamental to the program, constructed by Lusztig \cite{L90a}, for computing the characters of irreducible representations of algebraic groups over fields of positive characteristic  \cite{KL94}, \cite{AJS}, \cite{KT}. More recently, Kumar and Littelmann \cite{Li}, \cite{KL} succeeded in obtaining proofs of many theorems on the geometry of Schubert varieties, (including for example their normality) via the quantum Frobenius map and its splitting. Finally, the existence of the quantum Frobenius is also used in establishing the connection between quantum groups at a root of unity and perverse sheaves on the affine Grassmannian, \cite{ABG}.  

The original proof of the existence of $Fr$ was a tour de force computation with generators and relations \cite{L90}, reaching its more refined form in \cite[Chapter 35]{L93}. Unfortunately this does not give a conceptual reason for the existence of $Fr$. In this paper we attempt to remedy this by giving a new construction of the map $Fr$ using the language of Hall algebras. The Hall algebra construction realizes $\mathbf U^+$, the positive part of the quantum group, as an algebra of  functions on the moduli space of quiver representations. In this context, we show that $Fr_{|\mathbf U^+}$ can be interpreted as a restriction to the fixed-point set for the action of a non-split torus. Since the existence of $Fr$ on $\mathbf U_\varepsilon(\mathfrak g)$ follows easily from the existence of $Fr_{|\mathbf U^+}$ we readily recover the existence of $Fr$. In fact, following the approach of \cite{L93}, we use the modified form of the quantum group, as this form is better suited for specialization. The same strategy gives a geometric construction of the quantum Frobenius on $q$-Schur algebras following the ideas of \cite{BLM90}. Since the argument in this case is somewhat simpler we present it first. It also has the advantage of giving a construction of $Fr$ on all of $\mathbf U_\varepsilon(\mathfrak g)$ in the case $\mathfrak g = \mathfrak{sl}_n$. 

This work was motivated by a desire to gain a geometric understanding of the quantum Frobenius. As such, this paper constitutes the first step in that direction, since we work with functions over a finite field rather than \'etale sheaves, and the ``fasceaux-fonction'' lift remains to be done.

Section \ref{background} recalls the relevant definitions, and reduces the construction of the quantum Frobenius map to the positive part of the quantum group in the general case. Section \ref{qSchursection} gives the construction of $Fr$ in the case of $\mathfrak{sl}_n$ using $q$-Schur algebras. Section \ref{Hallalgebra} recalls the construction of the positive part of the quantum group using Hall algebras, and the final section concludes the construction of $Fr$. 

\textit{Acknowledgements}: The author was supported by a Royal Society University Research Fellowship while this paper was being written. 

\section{Quantum groups and the quantum Frobenius.}
\label{background}
In this section we recall the necessary background on quantum groups and the quantum Frobenius. Following \cite{L93} we will work with a modification of the quantum group which essentially replaces the ``Cartan'' part with a collection of orthogonal idempotents. Thus we begin with the definition of a quantum group and its associated modified form, following closely the presentation of \cite{L93}.
 
\begin{definition}
A \textit{Cartan datum} is a pair $(I,\cdot)$ consisting of a finite set
$I$ and a $\mathbb Z$-valued symmetric bilinear pairing on the
free Abelian group $\mathbb Z[I]$, such that
\begin{itemize}
\item $i\cdot i \in \{2,4,6,\ldots\}$
\item $2\frac{i\cdot j}{i\cdot i} \in \{0,-1,-2,\ldots\}$, for $i\neq j$.
\end{itemize}
We will write $a_{ij} = 2\frac{i\cdot j}{i \cdot i}$. Note that the matrix $A= (a_{ij})$ is a symmetrizable generalized Cartan matrix.
A \textit{root datum} of type $(I,\cdot)$ is a pair $Y,X$ of finitely-generated
 free Abelian groups and a perfect pairing
$\langle,\rangle \colon Y \times X \to \mathbb Z$, together with
imbeddings $I\subset X$, ($i\mapsto \alpha_i$) and $I\subset Y$, ($i
\mapsto \check{\alpha}_i$) such that $\langle \check{\alpha}_i,\alpha_j \rangle = 2\frac{i\cdot j}{i
\cdot i}$.  
\end{definition}

Let $v$ be an indeterminate. We will consider algebras over $\mathbb Q(v)$ and $\A = \mathbb Z[v,v^{-1}]$. For each $i \in I$ we set $v_i = v^{(i \cdot i)/2}$. Given only a Cartan datum, we may define an algebra $\mathbf f$ over $\mathbb Q(v)$ as follows (for more details see \cite[Chapter 3]{L93}): Take the free associative algebra $\mathscr F$ on generators $\{\theta_i: i \in I\}$. Then $\mathscr F$ is obviously $\mathbb Z[I]$-graded, and we denote the graded pieces by $\mathscr F_\nu,$ for $\nu \in \mathbb Z[I]$. For $x \in \mathscr F_\nu$ we write $|x|= \nu$. We may define an algebra structure on $\mathscr F \otimes \mathscr F$ by setting
\[
(x_1 \otimes x_2)(y_1\otimes y_2) =  v^{|x_2|\cdot|y_1|}(x_1y_1 \otimes x_2y_2).
\]
Let $r\colon \mathscr F \to \mathscr F \otimes \mathscr F$ be the homomorphism defined by setting 
\[
r(\theta_i) =  \theta_i \otimes 1+ 1\otimes \theta_i.
\]
It is straight-forward to show that there is a unique symmetric bilinear form $(\cdot, \cdot) \colon \mathscr F \times \mathscr F \to \mathbb Q(v)$ which satisfies
\begin{itemize}
\item $(\theta_i, \theta_j) = \delta_{ij} \frac{1}{1-v_i^2}$ for all $i,j \in I$.
\item $(x_1x_2, y) = (x_1 \otimes x_2, r(y))$ for all $x_1,x_2, y \in \mathscr F$.
\end{itemize}
(where the inner product on $\mathscr F \otimes \mathscr F$ is induced from that on $\mathscr F$ in the obvious way). Let $\mathscr I$ be the radical of this form, and define $\mathbf f$ to be the quotient $\mathscr F /\mathscr I$.

The quantum group $\mathbf U$ attached to a root datum $(X,Y)$ is the $\mathbb Q(v)$-algebra generated by symbols $E_i, F_i, K_\mu$, $i \in I$, $\mu \in Y$, subject to the following relations.
\begin{enumerate}
\item $K_0=1$, $K_{\mu_1}K_{\mu_2} = K_{\mu_1+\mu_2}$ for $\mu_1,\mu_2 \in Y$;
\item $K_{\mu} E_i K_{\mu}^{-1} = v^{\langle\mu,\alpha_i\rangle}E_i, \quad K_{\mu} F_i K_{\mu}^{-1} =
 v^{-\langle\mu,\alpha_i \rangle}F_i$ for all $i \in I$, $\mu \in Y$;
\item $E_iF_j - F_jE_i = \delta_{i,j}\frac{\tilde{K}_i-\tilde{K}_i^{-1}}{v_i-v_i^{-1}}$;
\item The maps $+ \colon \{\theta_i: i \in I\} \to \U$ given by $\theta_i \mapsto E_i$ and $-\colon \{\theta_i \in I\} \to \U$ given by $\theta_i \mapsto F_i$ extend to homomorphisms $\pm\colon \mathbf f \to \U$.

\end{enumerate}
Here $\tilde{K}_i$ denotes $K_{(i\cdot i/2)\check{\alpha}_i}$. The images of $\mathbf f$ under the maps $\pm$ are denoted $\mathbf U^\pm$. Note that this definition is the one used in \cite{L93} rather than the standard one involving the $q$-analogue of the Serre relations. However, the results of \cite[Chapter 33]{L93} show that these two definitions are equivalent.

For specialization to other coefficient rings, it is better to work with Lusztig's modified form of the algebra $\U$ denoted $\Ud$ in \cite{L93}. We briefly recall its construction. Let $\text{Mod}_X$ denote the category of left $\mathbf U$-modules endowed with a weight decomposition, thus the objects of $\text{Mod}_X$ are $\U$-modules $V$ such that
\[
V = \bigoplus_{\lambda \in X} V_\lambda,
\]
where
\[
V_\lambda = \{u \in V \colon K_\mu u = v^{\langle \mu,\lambda\rangle}u, \forall \mu \in Y\}.
\]
Let $\hat{\U}$ be the endomorphism ring of the forgetful functor from $\text{Mod}_X$ to the category of vector spaces. Thus by definition an element of $a$ of $\hat{\U}$ associates to each object $V$ of $\text{Mod}_X$ a linear map $a_V$, such that $a_W\circ f = f \circ a_V$ for any morphism $f\colon V \to W$.
Any element of $\mathbf U$ clearly determines an element of $\hat{\U}$, giving a natural inclusion $\U \hookrightarrow \hat{\U}$. For each $\lambda \in X$,
let $1_\lambda \in \hat{\U}$ be the projection to the $\lambda$ weight space. Then $\hat{\U}$ is isomorphic to the
direct product $\prod_{\lambda \in X}\mathbf U 1_\lambda$, and we set $\Ud$ to be the $\mathbb Q(v)$-subalgebra
$$\dot{\mathbf U} = \bigoplus_{\lambda \in X} \mathbf U 1_\lambda.$$
The algebra $\Ud$ does not have a multiplicative identity, but instead a collection $\{1_\lambda: \lambda \in X\}$ of orthogonal idempotents. It is clear that the category $\text{Mod}_X$ is equivalent to a category of modules for $\Ud$, the category of \textit{unital} modules. Thus provided the representations of $\mathbf U$ we study are weight modules, we may work with $\Ud$ or $\U$ interchangeably. 

There are two distinct integral forms for a quantum group which are usually considered -- with or without divided powers. In what follows we will need the integral form with divided powers. Set 
\[
[n]_i = (v_i^n - v_i^{-n})/(v_i-v_i^{-1}) \in \mathbb Z[v,v^{-1}],
\] 
and define 
\[
[n]_i! = [n]_i[n-1]_i.\ldots[1]_i; \qquad {n \brack k}_i = \frac{[n]_i!}{[k]_i! [n-k]_i!},
\]
(note that in fact ${n \brack k}_i \in \A$) and set $E_i^{(r)}$ to be the $v_i$-divided power $\frac{E_i^r}{[r]_i!}$ (with $E_i^{(0)} = 1$). By \cite{L93} the $\A$-subalgebra $\Ud_\A$ of $\Ud$ generated by $E_i^{(n)}1_\lambda, F_i^{(n)}1_\lambda$ $(i \in I, n \geq 0, \lambda \in X$) is an integral form (\textit{i.e.} the canonical map $\mathbb Q(v) \otimes_{\A} \Ud_\A \to \Ud$ is an isomorphism).  Similarly $\mathbf f_\mathcal A$, the $\mathcal A$-subalgebra of $\mathbf f$ generated by the elements $\theta_i^{(n)}$, ($i \in I, n \geq 0$), gives an $\mathcal A$-form of $\mathbf f$. Again by \cite{L93}, both $\Ud_\A$ and $\mathbf f_\mathcal A$ have canonical bases $\Bd$ and $\B$ respectively and hence in particular are free $\A$-modules (this is the main reason $\Ud$ is preferable to $\mathbf U$ when studying integral properties). For any $\A$-algebra $R$ we will write $_R\Ud$ for the specialization $R \otimes_{\mathcal A} \Ud_\mathcal A$ of $\Ud_\A$, and similarly write $_R\mathbf f$ and $_R\mathbf U^{\pm}$ for their corresponding specializations. The following lemma gives a presentation of $_R\Ud$. 

\begin{lemma}
\label{present}
For any $\A$-algebra $\phi\colon \A \to R$, the algebra $_R \Ud$ has the following presentation.  $_R \Ud$ is generated by elements $u^+1_\zeta u^-$ and $u^- 1_\zeta u^+$ for $u^\pm \in$ $_R\mathbf U^{\pm}$ and $\zeta \in X$, subject to the relations:
\begin{itemize}
\item $E_i^{(a)} 1_\zeta F_j^{(b)} = F_j^{(b)}1_{\zeta + a\alpha_i + b\alpha_j} E_i^{(a)}$ for $i \neq j$;

\item $E_i^{(a)} 1_{-\zeta} F_j^{(b)} = \sum_{t \geq 0} \phi({a+b-\langle \check{\alpha}_i, \zeta\rangle \brack t}_i) F_i^{(b-t)} 1_{-\zeta +(a+b-t)\alpha_i}E_i^{(a-t)}$;

\item $F_i^{(a)} 1_{\zeta} E_j^{(b)} = \sum_{t \geq 0} \phi({a+b-\langle \check{\alpha}_i, \zeta\rangle \brack t}_i) E_i^{(a-t)} 1_{\zeta -(a+b-t)\alpha_i}F_i^{(b-t)}$;

\item $(u^+ 1_\zeta)(1_{\zeta'}u^-) = \delta_{\zeta, \zeta'}u^+ 1_\zeta u^-$ and $(u^- 1_\zeta)(1_{\zeta'}u^+) = \delta_{\zeta, \zeta'}u^- 1_\zeta u^+$;

\item $_R \Ud$ is a left module for $_R \mathbf U^+$ and $_R \mathbf U^-$.
\end{itemize}
for $a, b \in \mathbb Z_{\geq 0}$, $\zeta,\zeta' \in X$, and $u^{\pm} \in \mathbf U^{\pm}$.
\end{lemma}
\begin{proof}
This is shown, in slightly different notation, in \cite[31.1.3]{L93}. The fact that the presentation holds for any ring $R$ is a consequence of the fact that $\mathbf f_{\mathcal A}$ is a free $\mathcal A$-module. 
\end{proof}

\begin{remark}
The fact that $\mathbf f_\A$ is a free $\A$-module is, to the author's knowledge, only known in complete generality via the existence of the canonical basis. In special cases (such as for finite type algebras) it can be shown by more elementary means.
\end{remark}

\begin{definition}
\label{cyclotomic}
Let $\ell$ be any positive integer.
We set $\A_\ell$ to be the quotient ring $\A/(\Phi_{2\ell}(v))$ where $\Phi_d$ is the $d$-th cyclotomic polynomial, and then let $\Ud_{\ell}$ and $\mathbf f_\ell$ be the corresponding specialization $\A_\ell \otimes_{\A} \Ud_{\A}$ and $\A_\ell \otimes_{\A} \mathbf f_{\A}$.
\end{definition}

Assume from now on that $\ell$ is relatively prime to the integers $\{\frac{1}{2}(i \cdot i): i \in I\}$. Let $X^\sharp = \{ \lambda \in X: \langle  \check{\alpha}_i, \lambda\rangle \in \ell\mathbb Z,  \forall i \in I\}$, and $Y^\sharp = \text{Hom}(X^\sharp, \mathbb Z)$. Define simple roots and coroots by setting $\alpha_i^\sharp = \ell\alpha_i$ and $\check{\alpha}_i^\sharp \in Y^\sharp$ to be the function whose value at $\lambda \in X^\sharp$ is $\ell^{-1}\langle \check{\alpha}_i, \lambda\rangle$. Then $(X^\sharp,Y^\sharp,I, \cdot)$ is a new root datum of type $(I, \cdot)$. Denote the associated modified quantum group by $\Ud^\sharp$. Note that  for any $\A$-algebra $R$, both $_R\Ud$ and $_R\Ud^\sharp$ are bimodules for $_R\mathbf f$ via both maps $\pm\colon {_R\mathbf f} \to {_R\mathbf U}^{\pm}$, since they are attached to the same Cartan datum. 

Let $\rho \colon \mathcal A \to \mathcal A$ be the map given by $v \mapsto (-1)^{\ell+1}v^\ell$. We will write $v_i^\sharp = \rho(v)^{(i \cdot i)/2}$, and ${a \brack b}_i^\sharp$ for the quantum binomial coefficient evaluated at $v_i^\sharp$ (rather than $v_i$). We now consider a number of specializations: let $\phi\colon \A \to \A_\ell$ denote the natural quotient map, and let $\mathbf f_\ell$ denote the specialization $\A_\ell \otimes_\phi \mathbf f_\A$, and $\Ud_\ell$ the specialization $\A_\ell \otimes_\phi \Ud_\A$. Next let $\varphi = \phi \circ \rho$, and let $\mathbf f_\ell^\sharp$ denote the specialization $\A \otimes_\varphi \mathbf f_\A$, and $\Ud^\sharp_\ell$ the specialization $\A\otimes_\varphi \Ud^\sharp_\A$. For clarity we will write $\theta_i^{(n)}$ and $E_i^{(n)}1_\lambda$ for the specializations of these elements in $\mathbf f_\ell$ and $\Ud_\ell$ respectively, while we will write $\vartheta_i^{(n)}$ and $e_i^{(n)}1_\lambda, f_i^{(n)}1_\lambda$ for their images in $\mathbf f_\ell^\sharp$ and $\Ud_\ell^\sharp$ respectively. The following theorem will be established in Section \ref{theproof}. The first version of such a theorem is due to  Lusztig \cite{L90}, while his most general version is given in \cite{L93}.

\begin{theorem}
\label{fFrob}
Suppose that $\ell$ is coprime to each of $\{\frac{1}{2}(i \cdot i): i \in I\}$.
There is a unique algebra homomorphism $Fr\colon \mathbf f_\ell \to \mathbf f_\ell^\sharp$ such that $Fr(\theta_i^{(n)}) = \vartheta_i^{(n/\ell)}$ if $n \in \ell \mathbb Z$ and $Fr(\theta_i^{(n)}) =0$ otherwise. 
\end{theorem}

The quantum Frobenius morphism is an analogous map between modified quantum groups, again the original construction of such a map was given by Lusztig \cite{L93}.
\begin{theorem}
\label{UFrob}
Let $R$ be an $\A$-algebra such that the homomorphism $\sigma\colon \A \to R$ factors through the map $\A \to \A_\ell$. Let $R_\rho$ denote $R$ viewed as an $\A$-algebra via the map $\sigma \circ \rho$. Then there is an unique surjective homomorphism $Fr\colon {_R\Ud} \to {_{R_\rho}\Ud^\sharp}$ such that 
\begin{enumerate}
\item $Fr(E_i^{(n)}1_\lambda) = e_i^{(n/\ell)}1_{\lambda}$ if $\lambda \in X^\sharp$ and $\ell$ divides $n$, and zero otherwise. 
\item $Fr(F_i^{(n)}1_\lambda) = f_i^{(n/\ell)}1_{\lambda}$ if $\lambda \in X^\sharp$ and $\ell$ divides $n$, and zero otherwise.
\end{enumerate}
\end{theorem}

We now show how Theorem \ref{UFrob} may be deduced from Theorem \ref{fFrob} and Lemma \ref{present}. For this we need some lemmas on Gaussian binomial coefficients. 

\begin{lemma}
\label{binom}
Let $l$ be a positive integer and $\psi\colon \A \to R$ be an algebra over $\A$ such that $\psi(v^{2l}) =1$, but $\psi(v^{2t}) \neq 1$ for all $0< t < l$. Then if $\mathbf v = \psi(v)$ we have 
\begin{enumerate}
\item $\psi({a \brack t})= 0$ if $l$ divides $a$ but not $t$.
\item Suppose that for $m \geq k$ we have $m = m_1l + s$ and $k = k_1l + t$ where $0 \leq s,t < l$.  Then we have
\[
\psi \big({m \brack k}\big)= \mathbf v^{l(k_1s - m_1 t)+(m_1+1)k_1l^2}{m_1 \choose k_1}\psi\big({s \brack t}\big).
\]

\end{enumerate}
\end{lemma}
\begin{proof}
See Chapter $34$ on Gaussian binomial coefficients in \cite{L93}. 
\end{proof}

\begin{lemma}
\label{checkingnormalization}
Let $\phi\colon \mathcal A \to \mathcal A_\ell$. Then if $a,t$ are divisible by $\ell$ we have $\phi( {a \brack t}_i) = \phi( {a/\ell \brack t/\ell}^\sharp_i)$, while if $a$ is divisible by $\ell$ and $t$ is not, we have $\phi({a \brack t}) = 0$.
\end{lemma}
\begin{proof}
Let $p_i\colon \A \to \A$ be defined by $p_i(v) = v_i$. Then from the previous lemma with $l= \ell$ and $\psi = \phi \circ  p_i$ we see that if $a = b\ell$ and $t = s\ell$ then $\phi({a \brack t}_i) = \mathbf v_i^{(b+1)s\ell^2}{b \choose s}$ where we set $\mathbf v_i = \phi(v_i)$. Applying the same lemma to the map $\phi\circ \rho \circ p_i$ and $\mathbf v_i^\sharp = \phi(\rho(v_i)) = \pm 1$, with $l = 1$ we get that $\phi \circ \rho({ b \brack s}_i) = \phi({b \brack s}_i^\sharp) = ((-1)^{\ell+1}\mathbf v^{\ell})^{\frac{1}{2}(i\cdot i)(b+1)s} {b \choose s}$ (where $\mathbf v = \phi(v)$). Thus it is enough to check that 
\[
\mathbf v_i^{(b+1)s\ell^2} =  ((-1)^{\ell+1}\mathbf v^{\ell})^{\frac{1}{2}(i\cdot i)(b+1)s}
\]
in $\mathcal A_\ell$. But this is clearly implied by the equation $(-1)^{\ell+1}\mathbf v^\ell = \mathbf v^{\ell^2}$ which holds in $\mathcal A_\ell$ since $\mathbf v^\ell = -1$. The case where $\ell$ divides $a$ but not $t$ is similar but easier. 
\end{proof}

\noindent
\textit{Proof of Theorem } \ref{UFrob}: Clearly it is enough to establish the case $R= \A_\ell$ with $\sigma$ the quotient map. The uniqueness and surjectivity follow from the description of $Fr$ on generators, so we need only to show its existence. Lemma \ref{present} describes $\Ud_\ell$ in terms of its $\pm$ bimodule structures for $\mathbf f_\ell$ and $\Ud^\sharp_\ell$ in terms of its bimodule structures over $\mathbf f_\ell^\sharp$. Thus assuming Theorem \ref{fFrob} it is enough to check the relations in Lemma \ref{present} are compatible with the map $Fr$. But this follows immediately from the Lemma \ref{checkingnormalization}. 
\qed

The bulk of the work in establishing the existence of the quantum Frobenius morphism then lies in the proof of Theorem \ref{fFrob} when $R = \mathcal A_\ell$. This purpose of this paper is to give an elementary construction of this map in the context of the Hall algebra. 

The original quantum Frobenius map \cite{L90} mapped to the Kostant-Chevalley form of the classical enveloping algebra, which is essentially equivalent to the specialization of $\Ud$ at $v =1$ (indeed in \cite{L08} Lusztig has shown how to construct the associated group scheme over $\mathbb Z$ from $\Ud$ at $v=1$ by taking an appropriate restricted dual). The algebra $\Ud^\sharp_\ell$ is very close to this algebra, since the parameters $v_i^\sharp$ are all equal to $\pm 1$ in $\mathcal A_\ell$. Lusztig calls this situation \textit{quasiclassical} and shows \cite{L93} that under mild hypotheses a quasiclassical specialization is in fact isomorphic to the classical ($v=1$) specialization. More precisely, we say that a Cartan datum is \textit{without odd cycles} if we cannot find a sequence $i_1, i_2, \ldots, i_p, i_{p+1} = i_1$ in $I$ such that $p \geq 3$ and $i_s \cdot i_{s+1} <0$ for all $s = 1,2, \ldots, p$.  We have the following theorem:

\begin{theorem}
\label{nooddcycles}
\cite[33.2]{L93}
Let $(Y,X)$ be a root datum such that the associated Cartan datum $(I, \cdot)$ has no odd cycles, and let $\phi\colon \mathcal A \to R$ be an $\mathcal A$-algebra such that $\phi(v_i) = \pm 1$ for all $i \in I$. Then if ${_{R_0} \Ud}$ denotes the specialization of $\Ud_\mathcal A$ obtained from the map $\varphi\colon \mathcal A \to R$ given by $\varphi(v) =1$, we have ${_{R_0}\Ud} \cong {_R \Ud}$.
\end{theorem}

\begin{remark}
The condition that $(I,\cdot)$ has no odd cycles is, for example, automatically satisfied for any Cartan datum of finite or affine type which does not have a component of affine type $A_n$ when $n$ is odd, however, in this case the datum is simply laced, and the isomorphism of the previous theorem may be checked directly (see also remark \ref{affinecase} at the end of the next section for this case). 
\end{remark}

We end this section by comparing the quantum Frobenius constructed here to the one in \cite[Chapter 35]{L93}. There, given a Cartan datum $(I, \cdot)$, a modification $(I, \circ)$ of the Cartan datum is used which, when $\ell$ is coprime to all the integers $\frac{1}{2}(i \cdot i)$, is related to $(I, \cdot)$ by $i \circ j = \ell^2 i \cdot j$. It is easy to see that $(X^\sharp, Y^\sharp, I, \circ)$ is a root datum of type $(I,\circ)$. Let $\mathbf f^*$ the algebra attached to the Cartan datum $(I, \circ)$ and $\Ud^*$ the modified quantum group attached to the root datum $(X^\sharp, Y^\sharp,I,\circ)$. Lusztig expresses his quantum Frobenius map as a surjective homomorphism from $\Ud_\ell$ to $\Ud^*_\ell =  A_\ell\otimes_{\A} \Ud^*_\A$.  For completeness we note the following identifications. Let $\mathbf f_\ell^* = \A_\ell \otimes_\A \mathbf f^*_\A$, and denote the generators of $\mathbf f^*$ by $(\theta_i^*)^{(n)}$.

\begin{lemma}
\label{oddcycles}
Assume that $(I, \cdot)$ is without odd cycles (see the remarks before the statement of Theorem \ref{nooddcycles}), and let $\theta_i^\sharp$ and $\theta_i^*$ be the generators of $\mathbf f^\sharp$ and $\mathbf f^*$ respectively. Then there is an isomorphism 
\[
\alpha \colon \mathbf f^\sharp_\ell \to  \mathbf f^*_\ell. 
\]
characterized by $\alpha(\theta_i^\sharp) = \theta_i^*$. Moreover, this induces a corresponding isomorphism between $\Ud_\ell^\sharp$ and $\Ud^*_\ell$.
\end{lemma}
\begin{proof}
Let $\phi\colon \mathcal A \to {\mathcal A}_\ell$ be quotient map as before. We write $v_i^* = v^{(i \circ i)/2}$ and $v_i^{\sharp} = \rho(v_i)$. We are in the quasiclassical case, since the parameters $\phi(v_i^{\sharp}) = \pm1$ and $\phi(v_i^*) = \pm 1$. As a consequence, if we assume that the Cartan datum has no odd cycles we may apply the results of \cite[33.2.2]{L93}, which show that the algebras $\mathbf f^\sharp_\ell$ and $\mathbf f^*_\ell$ are generated by the elements $\vartheta_i$ and $\theta_i^*$ respectively, and both have a presentation in terms of Serre relations. 

We may therefore define an isomorphism by specifying the action on the generators: the map $\vartheta_i \mapsto \theta_i^*$, ($i \in I$), extends to an isomorphism of algebras $\alpha\colon \mathbf f^\sharp_\ell \to \mathbf f^*_\ell$ by the presentation in terms of Serre relations, provided we check that $\phi(v_i^*) = \phi(v_i^\sharp)$. But this follows similarly to the proof of Lemma \ref{checkingnormalization}. Setting as before $\mathbf v = \phi(v)$, we have $\mathbf v^\ell = -1$ in $\mathcal A_\ell$, and thus:

\[
\phi(v_i^\sharp)= ((-1)^{\ell+1}\mathbf v^\ell)^{(i \cdot i)/2} = (-1)^{(\ell)(i \cdot i)/2}
\]
while 
\[
\phi(v_i^*) = \mathbf v^{(i \circ i)/2} = \mathbf v^{\ell^2(i \cdot i)/2} = (-1)^{\ell(i\cdot i)/2}
\]
as required. The construction of the corresponding isomorphism for the modified quantum groups $\Ud_\ell^*$ and $\Ud_\ell^\sharp$ the follows from the construction of $\alpha$ and Lemma \ref{present}, as in the deduction of Theorem \ref{UFrob}.
\end{proof}

In fact in \cite{L93} Lusztig's construction works more generally: there is a quantum Frobenius for $Fr \colon \mathbf U_\ell \to \mathbf U_\ell^*$ for an integer $\ell$ with very mild conditions on $\ell$ and the Cartan datum, although then the definition of the Cartan datum $(I,\circ)$ is more subtle. For $i \in I$, let $l_i$ be the smallest positive integer such that $l_i(i \cdot i) \in 2\ell\mathbb Z$. Then in general $(I,\circ)$ is defined by $i \circ j = l_il_j(i\cdot j)$. Lusztig's construction then requires that:
\begin{enumerate}
\item
$(I,\cdot)$ has no odd cycles;
\item 
for any $i \neq j$ in $I$ with $l_j \geq 2$ one has $l_i \geq -\langle \check{\alpha}_i, \alpha_j \rangle +1$.
\end{enumerate}
In the cases where the present paper applies, $l_i$ is equal to $\ell$ for all $i \in I$ (what we call the ``nondivisible case''), but neither condition is necessary for our argument (though we used the first condition to see that our map coincides with that of \cite{L93}). Thus, for example, the second condition fails for $G_2$ when $\ell =2$, so that our construction covers this case (which had already been checked directly by Lusztig) and a number of previously unknown cases, supporting Lusztig's hope \cite[\S 35.5.2]{L93} that the quantum Frobenius should exist with no restrictions. The divisible case where some $l_i \neq \ell$, which is established in \cite{L93} when the second of the above conditions holds seems much less clear in the context of the Hall algebra.

\section{The $q$-Schur algebra}
\label{qSchursection}

The goal of this section is to construct the quantum Frobenius homomorphim for $\Ud(\mathfrak{sl}_n)$. 
To do this we first construct geometrically the map induced by $Fr$ on $S_q(n,r)$, the ``$q$-Schur algebra'' which is a quotient of $\Ud(\mathfrak{sl}_n)$. The contents of this section give an alternative approach to some of the results of \cite{McG} in type $A$, but the proofs here are distinct from those in that paper, and indeed the results are more precise. We recall briefly the construction of $S_q(n,r)$. 

Let $V$ be an $r$-dimensional vector space over a field $\mathbf k$, and $n$ a positive integer. Let $\mathcal F^n$ denote the space of $n$-step partial flags in $V$, that is 
\[
\mathcal F^n = \{ (0=F_0 \subseteq F_1 \subseteq \ldots \subseteq F_n = V): F_i \text{ a subspace of } V\}.
\]
Then the group $GL(V)$ acts transitively on the components $\mathcal F^n$ which are indexed by the set $\mathcal C_{n,r} = \{(a_1,a_2, \ldots, a_n) \in \mathbb N^n: \sum_{i=1}^n a_i = r\}$. Moreover, $GL(V)$ also acts on $\mathcal F^n \times \mathcal F^n$ with finitely many orbits. Here the orbits are indexed by the set $\Theta_r$ of $n\times n$ matrices $(a_{ij})$ with nonnegative integer entries such that $\sum_{i,j} a_{ij} = r$. If $(F,F') \in \mathcal F^n \times \mathcal F^n$ then the orbit it lies in is indexed by the matrix $(a_{ij})$ where 
\[
a_{ij} = \dim\biggl(\frac{F_i\cap F_j'}{(F_{i-1}\cap
F_j')+(F_i\cap F_{j-1}')}\biggr).
\]
We write $\mathcal O_A$ for the orbit indexed by $A \in \Theta_r$. 

Now suppose that $\mathbf k = \mathbb F_q$, a finite field with $q$ elements. Let $S_{\mathbf k}(n,r)$ denote the set of $\mathbb Z$-valued $GL(V)$-invariant functions on $\mathcal F^n \times \mathcal F^n$. Then $S_{\mathbf k}(n,r)$ is an algebra under convolution, and moreover if we let $\mathbf 1 _A$ denote the indicator function for the $GL(V)$ orbit indexed by $A$ then $\{\mathbf 1_A: A \in \Theta_r\}$ is a $\mathbb Z$-basis of $S_{\mathbf k}(n,r)$ and the structure constants of $S_\mathbf k(n,r)$ with respect to this basis are polynomial in $q$. Hence we may define the $q$-Schur algebra $S_q(n,r)$ to be the $\mathbb Z[q]$-algebra with basis $\{1_A: A \in \Theta_r\}$ such that $S_{\mathbf k}(n,r)$ is the specialization of  $S_q(n,r)$ at  $q = |\mathbf k|$ for any finite field $\mathbf k$, where $1_A$ specializes to the indicator function $\mathbf 1_A$. To be consistent with previous section, we extend scalars from $\mathbb Z[q]$ to $\mathcal A$ by setting $q= v^2$. We will denote this extended algebra by $S_v(n,r)$. If $R$ is an $\A$-algebra, we let $S_R(n,r)$ denote the algebra $R \otimes_\A S_v(n,r)$. Let $(X,Y)$ be the root datum of type $SL_n$, and $\Ud$ the corresponding modified quantum group. Then the following is well known:

\begin{lemma}
$S_v(n,r)$ is a quotient of the algebra $\Ud$. 
\end{lemma}
\begin{proof}
This is the quantum analogue of Schur-Weyl duality. The proof in the context we describe is essentially in \cite{BLM90}.
\end{proof}

We now fix a positive integer $\ell$. Our construction of the quantum Frobenius hinges on the (trivial) observation that the $q$-Schur algebra construction works for \textit{all} finite fields -- and so in particular, it works both for the field $\mathbb F_q$ and its degree $\ell$ extension $\mathbb F_{q^\ell}$ -- and on the existence of the forgetful functor $A$ from the category of $\mathbb F_{q^\ell}$-vector spaces to the category of $\mathbb F_q$-vector spaces. Take $V$ to be an $r$-dimensional vector space over $\mathbb F_{q^\ell}$ and then let $W=A(V)$ be the $\mathbb F_q$-vector space obtained by forgetting the $\mathbb F_{q^\ell}$ structure on $V$. Then we can build the $\mathbb F_{q^\ell}$-Schur algebra $S_{\mathbb F_{q^\ell}}(n,r)$ on $V$ and the $\mathbb F_q$-Schur algebra $S_{\mathbb F_q}(n, \ell r)$ on $W$. There is an obvious inclusion
\[
\iota\colon \mathcal F^n_V \times \mathcal F^n_V \hookrightarrow \mathcal F^n_W \times \mathcal F^n_W
\]
from which we obtain a restriction map $\iota^* \colon S_{\mathbb F_q}(n, \ell r) \to S_{\mathbb F_{q^\ell}}(n,r)$. We claim that after specialization, this map is an algebra homomorphism. Pick $\e$ a square root of $q$, and extend scalars to the coefficient ring $\mathscr A_\ell = \mathbb Z[\e^{\pm1}]/(\Phi_{2\ell}(\e))$, where $\Phi_k$ denotes the $k$-th cyclotomic polynomial. We begin by noting the following simple property of cyclotomic polynomials.

\begin{lemma}
\label{oddevencyclotomic}
Let $\ell$ be a positive integer,  
\begin{enumerate}
\item if $\ell$ is even, then $\Phi_{\ell}(t^2) = \Phi_{2\ell}(t)$;
\item if $\ell$ is odd, then $\Phi_{\ell}(t^2) = \Phi_\ell(t)\Phi_{2\ell}(t)$.
\end{enumerate}
Thus if $s = \Phi_\ell(q)$, then for any $\ell$ we have $s=0$ in $\mathscr A_\ell$.
\end{lemma}
\begin{proof}
Both parts follow readily from the formula:
\[
\Phi_n(t) = \prod_{d | n} (x^d -1)^{\mu(n/d)},
\]
where $\mu$ is the Mobius function. 
\end{proof}

\begin{prop}
\label{qSchurFrob}
Let $I_r \colon \mathscr A_\ell \otimes_\mathbb Z S_{\mathbb F_q}(n,\ell r) \to \mathscr A_\ell \otimes_\mathbb Z S_{F_{q^\ell}}(n,r)$ be the map induced by $\iota^*$, then $I_r$ is an algebra homomorphism.
\end{prop}
\begin{proof}
We compute using the basis $\{\mathbf{1}_A: A \in \Theta_{\ell r}\}$. By simple dimension considerations, we see that 
\[
\iota^*(\mathbf{1}_D) = \left\{\begin{array}{cc}\mathbf{1}_{D'}, & \text{ if there is a } D' \in \Theta_r \text{ with } D= \ell D' \\0, & \text{ otherwise.}\end{array}\right.
\]
Suppose that $A, B \in \Theta_{\ell r}$, then we may write 
\[
\mathbf{1}_A . \mathbf{1}_B = \sum_{ C \in \Theta_{\ell r}} c_{A,B}^C \mathbf{1}_C,
\]
where if we fix $(F^0, F^1)$ in the orbit of $\mathcal O_C$, then $c_{A,B}^C$ is the number of points in the set
\[
S_{A,B}^C = \{ F \in \mathcal F_W \colon (F^0, F) \in \mathcal O_A, (F, F^1) \in \mathcal O_B\}.
\]
Hence applying $\iota^*$ to this equation, the terms on the right vanish unless $C = \ell C'$ for some $C' \in \Theta_r$, so that we get
\[
\iota^*(\mathbf{1}_{A}. \mathbf{1}_{B}) = \sum_{C' \in \Theta_r} c_{A, B}^{\ell C'} \mathbf{1}_{C'}.
\]
On the other hand, the product $\iota(\mathbf 1_A)\iota(\mathbf 1_B)$ is zero unless there are $A', B' \in \Theta_r$ such that $\ell A' = A$ and $\ell B' = B$, in which case it is $\sum_{C' \in \Theta_r} c_{A',B'}^{C'} \mathbf 1_{C'}$, where $\{c_{A',B'}^{C'}\}$ are the structure constants of $S_{\mathbb F_q^\ell}(n,r)$. 

By Lemma \ref{oddevencyclotomic} if $s = \Phi_\ell(q) \in \mathbb Z$, then $s=0$ in $\mathscr A_\ell$, hence to check that $\iota$ is an algebra homomorphism, it is enough to show that 

\begin{equation}
\label{congruence1}
c_{A, B}^{\ell C'} \equiv  \left\{\begin{array}{cc} c_{A',B'}^{C'} \mod s, & \text{ if } \exists A', B' \in \Theta_r \text{ with } A= \ell A', B = \ell B'; \\ 0 \mod s, & \text{ otherwise.}\end{array}\right.
\end{equation}
To show this we first let $T$ be the multiplicative group of $\mathbb F_{q^\ell}$. Then $T$ acts $\mathbb F_q$-linearly on $W$, and if $U$ is a subspace of $W$, it comes from an $\mathbb F_{q^\ell}$-subspace of $V$ if and only if it is preserved by the action of $T$. It follows from this that $T$ acts on $\mathcal F_W^n \times \mathcal F_W^n$ with fixed point set exactly equal to $\mathcal F^n_V \times \mathcal F^n_V$. Now when computing the coefficient $c_{A, B}^{\ell C'}$ we may assume that the flags $F^0$ and $F^1$ are fixed by $T$ (\textit{i.e.} that they are $\mathbb F_{q^\ell}$-subspaces), and hence the set $S_{A,B}^{\ell C'}$ has a $T$-action. The number of $T$-fixed points in $S_{A,B}^{\ell C'}$ is exactly the structure constant $c_{A',B'}^{C'}$ (or zero if no $A', B'$ exist). Hence using Equation \ref{congruence1} it is enough to show that 
\begin{equation}
\label{congruence}
|S_{A,B}^{\ell C'}| \equiv |(S_{A,B}^{\ell C'})^T| \mod s.
\end{equation}

We have already observed that  an $\mathbb F_q$-subspace of $W$ is an $\mathbb F_{q^{l}}$-subspace if and only if it is preserved by the multiplication action of the nonzero scalars $\mathbb F_{q^{l}}^*$ in $\mathbb F_{q^{l}}$, a cyclic group of order $q^{l} - 1$. Moreover, given an arbitrary $\mathbb F_q$-subspace $U$ of $W$, its stablizer in $T$ is exactly the multiplicative group of the largest subfield $\mathbb F_{q^d}$ of $\mathbb F_{q^l}$ which preserves it. Thus the size of the $\mathbb F_{q^l}^*$-orbit of $U$  is $(q^l-1)/(q^d -1)$ and so of order divisible by $s = \Phi_{l}(q)$ unless $U$ is a fixed point. Now \textit{a fortiori} this implies that the order of a nontrivial $\mathbb F_{q^\ell}$-orbit of $\mathbb F_q$-flags in $V$ is divisible by $s$ unless it is fixed, and hence Equation \ref{congruence} holds, proving the proposition.
\end{proof}

Since the proposition above holds for any prime power $q$, and so in particularly for infinitely many integers, it is then easy to see that the map $I_r$ must lift to a generic map. For an $\A$-algebra $\phi\colon \A \to R$, let $S_R^*(n,r)$ be the specialization of the $q$-Schur algebra corresponding to the map $\phi' \colon \A \to R$ given by $v \to \phi(v)^\ell$. 

\begin{cor}
There is a surjective homomorphism $Q_r \colon S_{\mathcal A_\ell}(n,r) \to S^*_{\mathcal A_\ell}(n,r)$.
\end{cor}
\begin{proof}
We define a map $Q_r \colon S_{\mathcal A_\ell}(n,\ell r) \to S^*_{\mathcal A_\ell}(n,r)$ by setting 
\[
Q_r(1_A) = \left\{\begin{array}{cc}1_{A'}, & \text{ if there is a } A' \in \Theta_r \text{ with } A= \ell A' \\0, & \text{ otherwise.}\end{array}\right.
\]
If we write the product
\[
1_A.1_B = \sum_{C \in \Theta_{\ell r}} c_{A,B}^C 1_C,
\]
for $c_{A,B}^C \in \A_\ell$, then, as above, it is clearly enough to show that
\[
1_{A'}1_{B'} = \sum_{D \in \Theta_{r}} c_{\ell A',\ell B'}^{\ell D} 1_D
\]
in $S^*_{\A_\ell}(n,r)$ (where terms which do not make sense are interpreted to be zero). Now by the previous proposition we know that this equation holds for infinitely many specializations of $\A_\ell$ (with $v^2 = q$ for any prime-power $q$) and so the result follows.
\end{proof}

It is easy to see that in fact $S^*_{\A_\ell}(n,r)$ is simply $\A_\ell\otimes_\mathbb Z S_\mathbb Z(n,r)$, thus in fact we have obtained a map to the classical Schur algebra (with scalars suitably extended).

\begin{remark}
\label{affinecase}
The algebras $S_q(n,r)$ can be put into an inverse system \cite{L99a}, the limit of which contains the algebra $\Ud$. The maps $Q_r$ are easily seen to be compatible with the maps in this inverse system. Thus, at least in type $A$ we get a construction of the entire map $Fr$ by ``geometric'' means. This reverses the logic of \cite{McG}, where it was shown, assuming its existence, that the quantum Frobenius induces a map between $q$-Schur algebras. However I do not know how to show the compatibility of $Fr$ with the basis $\{1_A: A \in \Theta_r\}$ by the methods of \cite{McG}. In that sense, the results of this section are more precise. Moreover, the method of this section can also equally be applied essentially word for word to the affine $q$-Schur algebras of \cite{L99}, \cite{L99a}, where one gets a construction of $Fr$ for the affine quantum group of type $\widehat{\mathfrak{sl}}_n$. (This includes the case where the Dynkin diagram is a single odd cycle, the only affine type where an odd cycle can occur.)
\end{remark}

\begin{remark}
The construction of $S_q(n,r)$ was a prototype for the later work by Nakajima on affine quantum groups using quiver varieties. It is natural to hope that some analog of the construction in this paper could be made in  equivariant $K$-homology to realize the quantum Frobenius in that context.  A direct path from the context of this paper and Nakajima's work is however far from straight-forward: Nakajima uses coherent sheaves on varieties which in type $A$ are the cotangent bundles of the varieties $\mathcal F_V\times \mathcal F_V$. These coherent sheaves are related to $\mathcal D$-modules on $\mathcal F_V\times \mathcal F_V$ via the machinery of mixed Hodge modules (see for example \cite{T}). In turn mixed Hodge modules are characteristic zero analogues of Weil sheaves on the corresponding varieties over the algebraic closure of $\mathbb F_p$, and one obtains functions on $\mathbb F_q$-points via the trace of Frobenius. Thus to connect the techniques of this paper to Nakajima's context via this path would require already as a first step a sheaf-theoretic version of our construction. On the other hand, motivated by this paper, in work with A. Oblomkov the author has studied by analogy another realization the quantum Frobenius on quantum affine $\mathfrak{sl}_2$ using coherent sheaves on $\mathbb P^1$.

\end{remark}

\section{The Hall algebra}
\label{Hallalgebra}

The purpose of this section is to review a construction of the algebra $\mathbf f$ using linear algebra over finite fields. This goes back to the seminal paper of Ringel \cite{R}, however our principal source will be the work of Lusztig \cite{L98}. 

We fix a prime $p$ and an algebraic closure $\mathbf k$ of the field $\mathbb F_p$. We also fix $\sqrt{p}$ a square root of $p$, and a finite subfield $\mathbb F_q$ of $\mathbf k$. We work with varieties defined over $\mathbb F_q$ and with integer-valued functions on their $\mathbb F_q$-rational points. Such functions have an elementary notion of pull-back and push-forward: indeed let $Z,W$ be arbitrary finite sets, and $f\colon Z \to \mathbb Z$ and $g \colon W \to \mathbb Z$ be integer-valued functions. For a map $\pi\colon Z \to W$, we set
\[
\pi_!(f)(x)= \sum_{z\colon \pi(z)=x} f(z),
\]
and 
\[
\pi^*(g)(z) = g(\pi(z)).
\]

Let $Q = (J,H,a)$ be a graph with vertex set $J$, edge set $H$ and an automorphism $a$, of order $d$. We pick an orientation $\Omega$ of our graph, in other words a pair of maps $s,t \colon H \to J$, where for $h \in H$ the set $\{s(h), t(h)\}$ is exactly the pair of vertices incident to $h$. We say that an orientation is admissible if it is compatible with $a$ in the sense that $s(a(h)) = a(s(h))$ and $t(a(h)) = a(t(h))$. In all cases we need to consider there is at least one admissible orientation. Let $I$ denote the set of $a$-orbits in $J$.

Let $\mathbb N J$ be the monoid of formal $\mathbb N$-linear combinations of the elements of $J$, and similarly let $\mathbb N I$ be the monoid of formal $\mathbb N$-linear combinations of the elements of $I$. It is convenient to identify $\mathbb N I$ with the submonoid of $\mathbb N J$ fixed by $a$. 

Let $\mathcal C'$ be the category of finite dimensional $J$-graded $\mathbf k$-vector spaces with an $\mathbb F_q$-structure (the morphisms are graded linear maps).  Let $F$ be the Frobenius morphism attached to the rational structure. For $V$ an object in $\mathcal C'$, we write the graded dimension
\[
|V| = \sum_{ i\in I} \dim(V_i) i \in \mathbb NJ.
\]
Let $\mathcal C$ be the category with objects $(V,a)$ where $V$ is an object of $\mathcal C'$ and $a\colon V\to V$ is an $\mathbb F_q$ linear map, such that $a(V_j) \subset V_{a(j)}$ and such that for any $k \in \mathbb N$ and any $j \in J$ such that $a^k(j) = j$ we have $a^k$ acts as the identity on $V_j$. For an object in $\mathcal C$ the dimension $|V|$ is an element of $\mathbb NI$. For each $\nu \in \mathbb NI$ we let $\mathcal C_\nu$ be the full subcategory whose objects have dimension $\nu$. 

For $V \in \mathcal C$ set $G_V = \prod_{j \in J} GL(V_j)$, and set 

\[
E_{V, \Omega} = \bigoplus_{h \in \Omega} \text{Hom}(V_{s(h)}, V_{t(h)}).
\]
We write $x = (x_h)_{h \in \Omega}$ for the elements of $E_V$. If there is no danger of confusion we will write $E_V$ instead of $E_{V,\Omega}$. The action of $a$ on $V$ induces a natural action of $a$ on $E_V$ given by insisting that the compositions

\xymatrix{& & & V_{s(h)} \ar[r]^{x_h} & V_{t(h)} \ar[r]^{a(t(h))} & V_{a(t(h))} }
\noindent and 

\xymatrix{& & & V_{s(h)} \ar[r]^{a(s(h))} & V_{t(h)} \ar[r]^{a(x)_{a(h)}} & V_{a(t(h))} }
\noindent coincide. This also induces an action of $a$ on $G_V$ which is compatible with the action of $G_V$ on $E_V$, that is $a(gx) = a(g)a(x)$ for all $g \in G_V$, $x \in E_V$. Let $F_a = a\circ F = F \circ a$. Notice that if $a$ has order $d$ as an automorphism of $Q$, then since $F_a^d = F^d \circ a^d = F^d$, the $F_a$-fixed points are all defined over $\mathbb F_{q^d}$, thus we can, if we wish, restrict our attention to $\mathbb F_{q^d}$-points.

Let $G_V^{F_a}$, $E_V^{F_a}$ \textit{etc.} denote the $\mathbb F_q$-rational points of these varieties with respect to the rational structure given by $F_a$. Let $\mathscr H_V$ be the Abelian group of all $G_V^{F_a}$-invariant functions on the spaces $E_V^{F_a}$. By $G_V$-invariance, if $V, V'$ are objects of $\mathcal C$ such that $|V| = |V'|$ the spaces $\mathscr H_V$ and $\mathscr H_{V'}$ are \textit{canonically} isomorphic, thus given $\nu \in \mathbb NI$ we may write $\mathscr H_\nu$ for the space $\mathscr H_V$ for any $V$ with $|V| = \nu$. Let $\mathscr H$ be the $\mathbb NI$- graded group
\[
\bigoplus_{ \nu \in \mathbb NI} \mathscr H_\nu.
\]
This direct sum has the structure of an algebra: Suppose that $\nu_1, \nu_2 \in \mathbb NI$, and that $f_1, f_2$ are  in $\mathscr H_{\nu_1},  \mathscr H_{\nu_2}$ respectively. Let $V$ be a graded vector space with $|V| = \nu_1 + \nu_2$, and let $W$ a subspace of dimension $\nu_1$, so that $T = V /W$ has dimension $\nu_2$. Let $P$ be the stabilizer of $W$ and $U$ the subgroup of $P$ whose elements induce the identity map on $W$ and $T$ (thus $P/U$ is isomorphic to $G_W \times G_T$). Let $K \subset E_V$ be subvariety of $x \in E_V$ such that $W$ is $x$-stable. Then consider the diagram

\xymatrix{& & & \bar{E}& \ar[l]_{p_1} E' \ar[r]^{p_2} & E'' \ar[r]^{p_3} &E_{V}.}
where $\bar{E}= E_W \times E_T$, $E' = G_V \times_U K$, $E'' = G_V \times_P K$ and the maps are the obvious ones. The diagram above is clearly compatible with $F_a$, hence we may take fixed points to obtain the diagram:

\xymatrix{& & & \bar{E}^{F_a}& \ar[l]_{p_1} E'^{F_a} \ar[r]^{p_2} & E''^{F_a}\ar[r]^{p_3} &E_{V}^{F_a}.}
\noindent 
where, using Lang's theorem, we see that
\[
\bar{E}^{F_a} = E_{W}^{F_a} \times E_T^{F_a}, \qquad E'^{F_a} = G_V^{F_a} \times_{P^{F_a}} K^{F_a}, \qquad E''^{F_a} = G_V^{F_a} \times_{U^{F_a}} K^{F_a}.
\]

Then since $p_2$ is a principal $P/U$-bundle, it is easy to see that if $f_i \in \mathscr H_{\nu_i}$ ($i = 1,2$) then there is a unique function $g$ on $E'$ such that $p_1^*(f_1 \boxtimes f_2) = p_2^*(g)$. We define 
\[
f_1\star f_2 = (p_3)_!(g).
\]
$\mathscr H$ is called the Hall algebra attached to the graph $(J,H,a, \Omega)$. 

In the case where the graph $(J,H,a)$ corresponds to a finite Dynkin diagram it is the entire Hall algebra that we are intereseted in, however the general case this algebra is much too large and we use instead a subalgebra (sometimes known as the composition algebra). 

\begin{definition}
Let $\mathcal X$ be the set of pairs $(\underline{i}, \underline{c})$ where $\underline{i} = (i_1, i_2, \ldots, i_m)$, is a sequence of elements of $I$, and $\underline{c} = (c_1,c_2, \ldots, c_m)$ is a sequence with $c_j \in \mathbb N$. A flag of type $(\underline{i}, \underline{c})$ is a filtration 
\[
\mathfrak f = (V= V^0 \supseteq V^1 \supseteq V^2 \supseteq \ldots \supseteq V^m=0)
\]
where $|V^{k-1}/V^{k}| = c_ki_k$. We let $\Pic$ be the variety of flags of type $(\underline{i}, \underline{c})$. Given $x \in E_V$ we say that $\mathfrak f$ is $x$-stable if $x(V^k) \subseteq V^k$. Finally we set 
\[
\bPic = \{(x, \mathfrak f)\in E_V \times \Phi_{\ic}: \mathfrak f \text{ is } x \text{-stable} \},
\]
and let $\pi_{\ic}\colon \bPic \to \Pic$ be the obvious map. Note that there is a natural action of $G_V$, the Frobenius $F$ and $a$ on all these varieties. The map $\pi_{\ic}$ is compatible with these actions, and so restricts to a map $\pi_{\ic} \colon \Psi_{\ic}^{F_a} \to \Phi_{\ic}^{F_a}$. Let $\mathbf 1_{\ic} = (\pi_{\ic})_!(1)$ where $1$ is the constant function on $\bPic^{F_a}$. Let $\mathcal F_V$ be the Abelian group of $G_V^{F_a}$-invariant functions on $E_V^{F_a}$ generated by the functions $\{\mathbf 1_{\ic}: \ic \in \mathcal X\}$. For each $V \in \mathcal C$, this is a finitely generated Abelian group (in fact, a finitely generated free Abelian group, since there is evidently no torsion). 
\end{definition}

Just as discussed above, we may write $\mathcal F_\nu$ instead of $\mathcal F_V$ where $|V| = \nu \in \mathbb NI$. It can be shown that
\[
\mathcal F = \bigoplus_{\nu \in \mathbb NI} \mathcal  F_\nu
\]
is a subalgebra of $(\mathscr H, \star)$. Indeed more precisely we have the following:

\begin{lemma}
\label{multiply}
Let $\icp, \icpp \in \mathcal X$, and suppose that  $\underline i' = (i_1',i_2',\ldots, i_{m'}')$, $\underline i'' = (i''_1, i''_2, \ldots, i''_{m''})$ and $\underline c' = (c_1',c_2', \ldots, c_{m'}')$, $\underline c'' = (c_1'', c_2'', \ldots, c_{m''}'')$. Set 
\[
\underline i = (i_1',i_2', \ldots, i_{m'}', i_1'', i_2'', \ldots, i_{m''}''), \underline c = (c_1',c_2', \ldots, c_{m'}', c_1'', c_2'', \ldots, c_{m''}'').
\]
Then we have
\[
\mathbf 1_{\icp} \star \mathbf 1_{\icpp} = \mathbf 1_{\ic}.
\]
\end{lemma}
\begin{proof}
This follows immediately from the definitions -- see \cite[1.11,1.18]{L98} for more details.
\end{proof}

The algebra structure we have defined depends on the orientation $\Omega$. To obtain an algebra which does not depend on the orientation, we must twist the multiplication by  a cocycle. The cost for doing this is that we must extend scalars from $\mathbb Z$ to $\mathscr A = \mathbb Z[q^{\frac{1}{2}}, q^{-\frac{1}{2}}]$ (of course, if $q$ is an even power of a prime, then this ring is again $\mathbb Z$). There are two choices for $\sqrt{q}$: either of  $\pm (\sqrt{p})^r$ where $q = p^r$. 
For $\nu, \mu \in \mathbb Z[I]$ we define 
\[
m(\nu, \mu) = \sum_{i \in I} \nu_i \mu_i + \sum_{h \in \Omega} \nu_{s(h)} \mu_{t(h)},
\]
and then define a new multiplication $\circ$ on $\mathcal F$ given by
\[
f_1 \circ f_2 = q^{-m(\nu_1, \nu_2)/2} f_1 \star f_2, \qquad (f_i \in \mathcal F_{\nu_i}, i=1,2)
\]
The algebra $(\mathcal F, \circ)$ is known as the twisted composition algebra. 

Given the datum of a quiver with automorphism we can construct a Cartan datum as follows: Let $I$ be the set of $a$-orbits in $J$. 
The pairing is given by

\begin{itemize}
\item $i\cdot j = - |\{h \in H: s(h) \in i \text{ and } t(h) \in j\}|$,  for $i \neq j$;
\item $i \cdot i = 2|\{v \in J: v \in i\}|$.
\end{itemize}

We have the following theorem of Lusztig\cite{L98}, a generalization of the work of Ringel, which relates the twisted composition algebra to quantum groups. Let $\mathbf f$ be the algebra attached to the Cartan datum above in the manner of Section \ref{background}.

\begin{theorem}\cite[Theorem 1.20]{L98}
\label{Hall}
The algebra $(\mathcal F, \circ)$ is isomorphic to $\mathbf f_{\A |v = \sqrt{q}}$.
\end{theorem}
\begin{proof}
Theorem 1.20 in \cite{L98} actually deals with the algebras defined over the field $\bar{\mathbb Q}_l$. In order to obtain the statement over $\A$ one must combine that result with the results on the in Section 2 of that paper, where it is shown how to reconstruct to the algebra $\mathbf f_\A$ from the family of specializations obtained via the Hall algebra construction.
\end{proof}

\begin{remark}
In the finite type case, the algebra $\mathcal F$ is equal to the full Hall algebra $\mathscr H$, and this theorem is essentially a proof of the existence of Hall polynomials. In general, the full Hall algebra depends on the field $\mathbb F_q$ more finely, since the isomorphism classes of representations depend on the size of the field thus the notion of Hall polynomials does not automatically make sense. The above theorem shows in a quite general context that something like Hall polynomials make sense for the composition subalgebra. 
\end{remark}

In this paper, the algebras $\mathbf f$ are the objects of primary interest, thus it should be pointed out that given a Cartan datum $(I, \cdot)$, one can construct a quiver with automorphism whose associated Cartan datum is $(I, \cdot)$. Such a construction is not unique, so for definiteness we give a procedure for constructing such a graph from $(I, \cdot)$, following \cite[14.1]{L93}. Let $d_i = \frac{1}{2}(i \cdot i)$ and let $d$ be the least common multiple of the numbers $\{d_i: i \in I\}$, the lacing number of the Cartan datum. For each $i \in I$ let $D_i$ be a set with $d_i$ elements equipped with a cyclic action $a \colon D_i \to D_i$. For each pair $i,j$ with $i \cdot j <0$ consider the action of $a\times a $ on $D_i \times D_j$. Each $a\times a$-orbit $\mathcal O$ has size $\text{l.c.m.}(d_i,d_j)$ the least common multiple of $d_i$ and $d_j$. By definition, this divides $-i\cdot j$. Thus we can construct a set $H_{ij}$ of $(-i \cdot j)/\text{l.c.m.}(d_i,d_j)$ copies of $\mathcal O$ with a permutation $a$ acting as $a\times a$ on each copy of $\mathcal O$. 
We have a natural map from $H_{ij}$ to $D_i \times D_j$, and so setting $J= \sqcup_{i \in I} D_i$, $H = \sqcup_{i,j \in I, i\cdot j <0} H_{ij}$ we have a graph $(J,H,a)$ with automorphism with the required properties.

Finally, we give an explicit description of the structure of the $\mathbb F_q$-rational points of $E_V$ in the case of a nontrivial automorphism (when $a$ is trivial, this is clear). This is essentially equivalent to relating Lusztig's construction with the previous work of Dlab and Ringel on species \cite{DR}. (In the finite type case this is already given by Lusztig in \cite[\S 11]{L90b}). For each $i \in I$ pick $j_i \in J$ such that $j_i \in i$, and similarly for each $a$-orbit $k$ in $\Omega$ pick a representative $h_k$ (in the situation of the above construction, one can clearly arrange the choice of $h_k$ so that $s(h_k)$ and $t(h_k)$ are the representatives of their respective $a$-orbits). We have
\[
V^{F_a} = \bigoplus_{i \in I} V_i^{F_a},
\] 
where $V_i = \bigoplus_{j \in i} V_j$, since the $V_i$ are $F_a$-stable. Moreover, each $V_i^{F_a}$ naturally has the structure of a $\mathbb F_{q^{d_i}}$-vector space. To see this, observe that $V_i^{F_a}$ lies in $V_i^{F^{d_i}}$, and the projection onto any of the factors in the decomposition $V_i^{F^{d_i}} = \bigoplus_{j \in i} V_j^{F^{d_i}}$ gives a bijection between $V_i^{F_a}$ and $V_j^{F^{d_i}}$, which intertwines the action of $a^{-1}$ with the action of $F$.

The elements of $E_V^{F_a}$ consist of $\mathbb F_q$-linear maps $(x_k\colon V^{F_a}_{s(k)} \to V^{F_a}_{t(k)})$ where $k$ runs over the set of $a$-orbits in $\Omega$ and $s(k)$, $t(k)$ are the corresponding $a$-orbits in $J$. Suppose that $i_1$ and $i_2$ are $a$-orbits in $J$ with representatives $j_1$ and $j_2$ respectively and $k$ is an $a$-orbit of edges between the orbits $i_1$ and $i_2$ with representative $h$. Then $x_k$ corresponds to 
a map from 
\[
V_{j_1}^{F^{d_{i_1}}} \otimes_{\mathbb F_{q^{d_{i_1}}}} \mathbb F_{q^{- i_1\cdot i_2}} \to V_{j_2}^{F^{d_{i_2}}}
\]
where $\mathbb F_{q^{-i_1 \cdot i_2}}$ is a $\mathbb F_{q^{d_{i_1}}}$-$\mathbb F_{q^{d_{i_2}}}$ bimodule via 
\[
(\lambda, \mu)t = \lambda^{q^{s_1}} \mu^{q^{s_2}} t
\]
where $s(h) = a^{s_2}(j_1)$ and $t(h) = a^{s_1}(j_2)$. In fact, by the parenthetical remark, in the case of our explicit construction of a quiver attached to a Cartan datum, we may assume that $s_1=s_2=0$, so that the bimodule structure of $\mathbb F_{q^{-i_1\cdot i_2}}$ is the natural one. For the rest of the paper, we will assume that this is the case. (For a more detailed exposition of this translation, see, for example, \cite{DD}.)

\section{Quantum Frobenius}
\label{theproof}
The goal of this section is to construct the quantum Frobenius homomorphim in the context of the Hall algebra. Let $(I,\cdot)$ be a Cartan datum, and $\mathbf f$ the associated algebra. Fix a positive integer $\ell$, and assume that the prime $p$ of the previous section does not divide $\ell$ and that $\ell$ is relatively prime to the integers $\{\frac{1}{2}(i\cdot i): i \in I\}$. We will denote by $\mathbf f$ the algebra attached to $(I,\cdot)$ and $\mathbf f_\ell$ the algebra attached to $(I, \circ)$, the $\ell$-modified Cartan datum. Let $(J,H,a)$ be a quiver with an automorphism, such that its associated Cartan datum is $(I, \cdot)$. For convenience, we write $d_i = (i \cdot i)/2$. Let $\mathcal C$ be the category of $J$-graded $\mathbf k$-vector spaces with an $\mathbb F_q$-rational structure and an action of $a$, and $\mathcal C^\ell$ the corresponding category equipped with an $\mathbb F_{q^\ell}$ structure instead of an $\mathbb F_q$-structure.  Since our argument is more naturally phrased in the explicit description of the rational points given at the end of Section \ref{Hallalgebra}, we will from now on use that context. We will therefore also now, by abuse of notation, write $E_V$, $G_V$ etc. to mean the set of rational points.

The restriction of scalars functor from $\mathbb F_{q^\ell}$-vector spaces to $\mathbb F_q$-vector spaces extends to gives a natural functor $R_\ell$ from $\mathcal C^\ell$ to $\mathcal C$: namely each $\mathbb F_{q^{\ell d_i}}$-vector space corresponding to $i \in I$ ``forgets'' to yield an $\mathbb F_{q^{d_i}}$-vector space. The key point, once automorphisms are introduced, is that the field $\mathbb F_{q^{\ell d_i}}$ is simple as an $\mathbb F_{q^\ell}$-$\mathbb F_{q^{d_i}}$ bimodule when $\ell$ and $d_i$ are coprime.
 
If $V$ is an object in $\mathcal C^\ell$, then set $W = R_\ell(V)$. Clearly, an element $x \in E_V$ induces a map on $W$, thus we obtain a map
\[
\iota_\ell \colon E_V \to E_{W}.
\]

As in the case of the $q$-Schur algebra, we will need another description of the image of $E_V$ in $E_{W}$.  Let $T$ denote the multiplicative group of $\mathbb F_{q^{\ell}}$.  Then $T$ acts on $V$ by the obvious ``diagonal'' embedding of $\mathbb F_{q^\ell}^*$ in $G_V$. Then $T$ still acts on $W = R(V)$, and hence on $E_W$. Clearly we have the following lemma:

\begin{lemma}
\label{fixedpoints}
The image of $\iota_\ell$ is precisely the fixed-point set of the action of $T$, that is, 
\[
\iota_\ell(E_V) = (E_{R(V)})^T
\]
\end{lemma}
\begin{proof}
This follows immediately from the fact that an $\mathbb F_q$-linear map is $\mathbb F_{q^{\ell s}}$-linear if it is $\mathbb F_{q^{s}}$-linearly and commutes with the action of $T=\mathbb F_{q^\ell}^\times$, when $\ell$ and $s$ are coprime.
\end{proof}

 Let $\mathscr H$ and $\mathscr H^\sharp$ be the Hall algebras constructed using the categories $\mathcal C$ and $\mathcal C^\ell$ respectively, with twisted composition algebras $\mathcal F$ and $\mathcal F^\sharp$ respectively, where we pick $\e$ a square root of $q$ and pick $\e_\ell = (-1)^{\ell+1}\e^\ell$ as our square root of $q^\ell$. Thus $\mathcal F^\sharp$ and $\mathcal F$ are specializations of the algebra $\mathbf f_\mathcal A$ at $v = \e$ and $\e_\ell$ respectively. 

\begin{definition}
Let $Q \colon \mathscr H \to \mathscr H^\sharp$ be the $\mathscr A$-linear map given as follows: Let $f \in \mathcal F_\nu$. If $\nu \notin \ell\mathbb NI$, then set  $Q(f) = 0$. Otherwise we may choose a $V$ in $\mathcal C^\sharp$ such that $R_\ell(V)$ has dimension $\nu$, and we have the embedding $\iota_\ell \colon E_V \to E_{R_\ell(V)}$.
Set $Q(f) = \iota_\ell^*(f)$.
\end{definition}

Although the construction of $R_\ell$ involves choices, the map $Q$ on the Hall algebra is clearly independent of them, since the functions in $\mathcal H_V$ are $G_V$-invariant.
It is important to note that $Q$ is \textit{not} an algebra homomorphism, and indeed it is not even clear that it restricts to give a linear map between the corresponding composition subalgebras. In order to obtain these properties we must specialize the ring $\mathscr A$. Recall from Section \ref{background} that $\Phi_\ell$ is the $\ell$-the cyclotomic polynomial. Let $s = \Phi_\ell(q)$, and, as in the $q$-Schur algebra case, set $\mathscr A_\ell$ to be the ring $\mathbb Z[\e, \e^{-1}]/(\Phi_{2\ell}(\e))$.

In order to prove that $Q$ is an algebra homomorphism after we pass to $\mathscr A_\ell$, we use the generators $\{\mathbf 1_{\ic}: \ic \in \mathcal X\}$. However, because we have twisted the multiplication, it is convenient to renormalize them: 
For $\ic \in \mathcal X$, set 

\begin{equation}
\label{Nformula}
\begin{split}
N\ic & = \sum_{h \in \Omega; r < r': s(h) \in i_l, t(h) \in i_{l'}} c_r c_{r'} + \sum_{ j; r < r': j \in i_r = i_{r'} } c_r c_{r'} \\
& = \sum_{ r < r': i_r \cdot i_{r'} < 0} (i_r \cdot i_{r'})c_r c_{r'} + \sum_{ r < r': i_r = i_{r'} } (\frac{i_r \cdot i_r}{2})c_r c_{r'} \\
\end{split}
\end{equation}
and let $\vartheta_{\ic} = \e^{-N\ic}\mathbf 1_{\ic}$ (and similarly $\vartheta^\sharp_{\ic} = (\e_{\ell})^{- N\ic}\mathbf 1_{\ic}$).

\begin{lemma}
\label{residue}
Suppose that $\ell$ is an odd integer. Let $\ic \in \mathcal X$ with $\underline i = (i_1, i_2, \ldots, i_m)$ and $\underline c = (c_1, c_2, \ldots, c_m)$, and suppose $c_k = \ell b_k$ for all $k$, $ 1 \leq k \leq m$. We have
\[
Q(\vartheta_{\ic}) =  \vartheta_{(\underline i, \underline b)}^\sharp \in \mathscr A_\ell,
\]
Moreover, if $c_k \notin \ell \mathbb N$ for any $k$, we have
\[
Q(\vartheta_{\ic}) = 0 \in \mathscr A_\ell.
\]
Thus over $\mathscr A_\ell$ the map $Q$ restricts to a map $Q\colon \mathcal F \to \mathcal F^\sharp$.
\end{lemma}

\begin{proof}
We may assume that, for each $i \in I$ we have $\sum_{k: i_k =i}c_k \in \ell\mathbb N$ since otherwise it is clear that $Q(\mathbf 1_{\ic})$ vanishes. Thus we can find a vector space $V \in \mathcal C^\ell$ such that $W = R_\ell(V)$ has dimension $|W| = \sum_{k=1}^m c_k i_k$. We prove the lemma using the action of the group $T = \mathbb F_{q^{\ell}}^\times$ on $W$. The value of $\mathbf 1_{\ic}$ at a point $x \in (E_W)^T = \iota_\ell(E_V)$ is by definition the number of points in the set 
\[
S_x = \{\mathfrak f \in \Phi_{\ic} : \mathfrak f \text{ is } x\text{-stable}\}.
\]
Since $x$ is $T$-stable, this set has an action of $T$.
An $I$-graded subspace of $W$ is stable under $T$ if and only if it comes from a subspace of $V$, by the assumption that $\ell$ is coprime to each $d_i, i \in I$ (\textit{c.f.} the proof of Lemma \ref{fixedpoints}). Hence the value of $1_{(\underline i, \underline b)}$ at $x$ is $|S_x^T|$. 

Consider the $T$-orbits on $S_x$. As in the proof of Proposition \ref{qSchurFrob}, the nontrivial orbits all have stabilizers of order $q^{d} -1$ where $d$ is a proper divsior of $\ell$. Since for all proper divisors $r$ of $\ell$ we have that $s = \Phi_\ell(q)$ divides $(q^\ell-1)/(q^r-1)$, it follows that the order of the orbit of $U$ is divisible by $s$ unless $U$ is fixed by $T$. Since the stabilizer of a flag is the intersection of the stabilizers of the subspaces it contains, it is immediate that 
\[
|S_x| \equiv |S_x^T| \mod (s).
\]
Now Lemma \ref{oddevencyclotomic} shows that $Q$ is compatible with the untwisted algebra structure on the composition algebras. 

To show that $Q$ is compatible with the twisted multiplication, observe that by Equation (\ref{Nformula}) we have $N\ic = \ell^2N(\underline i, \underline b)$. But in $\mathscr A_\ell$ we have $\e^\ell = -1$, so that $\e^{-\ell^2N(\underline i, \underline b)} = \e_\ell^{-N(\underline i, \underline b)}$, if $\e^{\ell^2} = (-1)^{\ell+1}\e^{\ell}$ which is clear (both sides are $-1$ if $\ell$ is odd and $1$ if $\ell$ is even). 
\end{proof}

\begin{remark}
The last paragraph of the proof is where we needed to be careful in our choice of square roots.
\end{remark}

\begin{theorem}
\label{homo}
The map $Q\colon \mathcal F \to \mathcal F^\sharp$ is an algebra homomorphism over the ring $\mathscr A_\ell$. 
\end{theorem}
\begin{proof}
We must show that the multiplication is preserved. Since the algebra $\mathcal F$ is spanned as an $\mathscr A_\ell$-module by the elements $\theta_{\ic}$ where $\ic \in \mathcal X$ it is enough to show that 
\begin{equation}
\label{multiplication}
Q(\vartheta_{\icp} \circ \vartheta_{\icpp}) = Q(\vartheta_{\icp})\circ Q(\vartheta_{\icpp}).
\end{equation}
for any $\icp, \icpp \in \mathcal X$. The analog of Lemma \ref{multiply} for the twisted multiplication shows that $\vartheta_{\icp} \circ \vartheta_{\icpp} = \vartheta_{\ic}$ (where $\ic$ is as in the statement of that lemma), and so also that $\vartheta_{\icp}^\sharp \circ \vartheta_{\icpp}^\sharp = \vartheta_{\ic}^\sharp$. Equation \ref{multiplication} now follows immediately from Lemma \ref{residue}.
 \end{proof}

We may now establish Theorem \ref{fFrob} from Section \ref{background}. As in that section, we let $\rho \colon \mathcal A \to \mathcal A_\ell$ be the map given by $v \mapsto (-1)^{\ell+1}v^{\ell} \in \mathcal A_\ell$ and let $\mathbf f^\sharp_\ell = \mathcal A_\ell \otimes_{\rho} \mathbf f$. Thus the algebra $\mathcal F^\sharp$ is a specialization of $\mathbf f_\ell^\sharp$.

We start with a simple lemma.
\begin{lemma}
\label{vanishing}
Let $\mathcal S\subseteq \mathbb N$ be an infinite set of positive integers, and let $\mathcal A_\ell$ be the ring of cyclotomic integers $\mathbb Z[\zeta]$ where $\zeta = e^{\pi i/\ell}$ is a primitive $2\ell$-th root of unity. Then we have
\[
\bigcap_{q \in \mathcal S} (\zeta^2-q) = 0,
\]
where for $a \in \mathcal A_\ell$, we write $(a)$ for the principal ideal generated by $a$.
\end{lemma}
\begin{proof}
Suppose for the sake of contradiction that there is some $a \neq 0$ which lies in all of the ideals $(\zeta^2 -q)$, $(q \in \mathcal S)$. Let $N$ denote the norm for the field extension $\mathbb Q(\zeta) \supset \mathbb Q$. Since the norm is multiplicative, for any $q \in \mathcal S$ the norm $N(\zeta^2 - q)$ divides $N(a)$, and so in particular the set $\{N(\zeta^2-q): q \in \mathcal S\}$ is bounded. 

On the other hand, let $G$ be the Galois group of the extension $\mathbb Q(\zeta) \supset \mathbb Q$, so that $G \cong (\mathbb Z/\ell \mathbb Z)^\times$. Then by definition we have
\[
N(\zeta^2-q) = \prod_{g\in G}(g(\zeta)^2-q) = (-1)^{|G|}\prod_{g \in G}(q-g(\zeta)^2)
\]
If $\ell$ is odd, then $\zeta^2$ is a primitive $2\ell$-th root of unity, so that $N(\zeta^2-q) = \pm\Phi_{2\ell}(q)$, while if $\ell$ is even, then $\zeta^2$ is a primitive $\ell$-th root of unity, and so $N(\zeta^2-q) = \Phi_\ell(q)^2$.
Since $\Phi_\ell(t)$ and $\Phi_{2\ell}(t)$ both tend to infinity as $t$ does, and the set $\mathcal S$ must be unbounded since it is infinite, we have a contradiction.
\end{proof}

\begin{theorem}
\label{redFrob}
There is a unique algebra homomorphism $Fr\colon \mathbf f_\ell \to \mathbf f_\ell^\sharp$ such that $Fr(\theta_i^{(n)}) = \vartheta_i^{(n/\ell)}$ if $n \in l_i \mathbb Z$ and zero otherwise. 
\end{theorem}
\begin{proof}
We first define a map $d_\ell \colon \mathcal X \to \mathcal X\cup \{\emptyset\}$ by setting 
\[
d_{\ell}(\ic) = \left\{\begin{array}{cc}(\underline i, \underline b), & \text{if } c_k = l_{i_k}b_k, b_k \in \mathbb N \text{ for each } k \\\emptyset, & \text{otherwise}. \end{array}\right.
\]
It is known that the algebra $\mathbf f_\A$ has a canonical basis $\mathbf B$, and thus its specialization gives a basis of $\mathbf f_\ell$, which we denote by $\mathbf B_\ell$. Similarly $\mathbf B$ yields a basis $\mathbf B^\sharp_\ell$ of $\mathbf f_\ell^\sharp$. 
Since the monomials $\{\theta_{\ic}: \ic \in \mathcal X\}$ span $\mathbf f_\ell$, given a basis element $b \in \mathbf B_\ell$ we may write $b = \sum_{\ic \in K} c_{\ic} \theta_{\ic}$ for some finite set $K \in \mathcal X$ and coefficients $c_{\ic} \in \A$. We set 
\[
Fr(b) = \sum_{\ic \in K} c_{\ic} \vartheta_{d_\ell\ic},
\]
where we define $\vartheta_{\emptyset} = 0$.
Doing this for all $b \in \mathbf B$ and extending linearly we obtain a map $Fr\colon \mathbf f_\ell \to \mathbf f_\ell^\sharp$. We claim that this map is an algebra homomorphism. Suppose that $b,b' \in \mathbf B$ and consider $Fr(bb') - Fr(b)Fr(b')$. We may write this as a linear combination of the basis elements $b^\sharp \in \mathbf B^\sharp_\ell$, say
\[
Fr(bb') - Fr(b)Fr(b') = \sum_{b^\sharp \in \mathbf B^\sharp} c_{b^\sharp} b^\sharp
\]
for some $c_{b^\sharp} \in \A_\ell$ (all but finitely many of which are equal to zero). Thus we must show that the coefficients $c_{b^\sharp}$ are all zero. Let $q$ be a prime power. We may specialize $\mathbf f_\ell^\sharp$ to the ring $\Al / (v^2-q)\Al = \mathscr A_\ell$ (since both rings are just $\A/(\Phi_{2\ell}(v), v^2-q)$). Using Theorem \ref{Hall} and the definition of $Fr$ it follows that $Fr(b)= Q(b)$, and so since $Q$ is an algebra homomorphism we see that $c_{b^*}$ is an element of the ideal $(v^2-q) \subset \mathcal A_\ell$. But since this last condition must hold for infinitely many prime powers $q$ it follows from Lemma \ref{vanishing} that $c_{b^*} = 0$ as required. It is immediate that $Fr$ is in fact unique, since it is a homomorphism, and its value on the generators $\theta_i^{(n)}$ is uniquely determined. 
\end{proof}

\begin{remark}
The previous proof uses the existence of the canonical basis, but in a rather weak sense. Only the fact that  $\mathbf f_\mathcal A$ is a free $\mathcal A$-module is necessary for the argument. As mentioned before, for special cases, \textit{e.g.} in finite type, this is possible to establish without the canonical basis. 
\end{remark}

\begin{remark}
It is also known that $Fr$ respects the coproduct. Though this is also straightforward to check directly on the generators, since the twisted coproduct can be interpreted in the context of the Hall algebra, one can show directly its compatibility with $Q$. 

In his work on $Fr$, Lusztig also constructs a map $Fr'$ which acts as a kind of splitting map for $Fr$. From the point of view of this paper, that map is much more mysterious than the map $Fr$. 
\end{remark}

\begin{remark}

Lusztig's construction of canonical bases \cite{L90b} arises naturally from his lifting of Ringel's construction of $\mathbf U^+$ to the level of sheaves on the moduli space of representations of a quiver. The construction of this paper is at the level of Hall algebras, but if it could also be lifted to the level of perverse sheaves, it would give a context in which to study the compatibility of (the specializations of) the canonical basis with the map $Fr$. In low rank examples where the canonical basis has been explicitly computed, there is good evidence for such compatibility: indeed it can be checked that $Fr$ is compatible with the canonical basis, in the sense that a basis element is either mapped to zero or a basis element.

\end{remark}

\end{document}